\documentclass{tac} \relax

\DeclareFontFamily{U}{msa}{} \DeclareFontShape{U}{msa}{m}{n}{<-6>msam5<6
-8.5>msam7<8.5->msam10}{} \DeclareSymbolFont{latex-font
msa}{U}{msa}{m}{n} \DeclareFontFamily{U}{msb}{}
\DeclareFontShape{U}{msb}{m}{n}{<-6>msbm5<6-8.5>msbm7<8.5->msbm10}{}
\DeclareSymbolFont{latex-font msb}{U}{msb}{m}{n}

\DeclareMathAlphabet{\mathscr}{U}{rsfs}{m}{n}

\DeclareFontFamily{U}{cmr}{} \DeclareFontShape{U}{cmr}{a}{n}{<6>cmb6<8>c
mb8<12>cmb12}{} \DeclareMathAlphabet{\mathbt}{U}{cmr}{a}{n}

\DeclareMathSymbol\emptyset\mathord{latex-font msb}{"3F}
\DeclareMathSymbol\ZZ\mathord{latex-font msb}{`Z}

\csname@@input\endcsname xypic \space\newdir{ ^(}{{}*!/-4pt/\dir^{(}}

\title{On reflective subcategories of locally presentable categories}
\author{J.\ Ad\'amek and J.\ Rosick\'y\footnote{\supportedby}}
\def\supportedby{Supported by the grant agency of the Czech Republic
under reference P201/12/G028.} 
\dedication{Dedicated to the memory of Horst Herrlich} 
\address{J.\ Ad\'amek\\Institut f\"ur Theoretische
Informatik\\Technische Universit\"at Braunschweig\\38032 Braunschweig,
Germany\\\texttt{j.adamek@tu-bs.de}\\\\J.\ Rosick\'y\\Department of
Mathematics and Statistics\\Masaryk University, Faculty of
Sciences\\Kotl\'a\v{r}sk\'a 2, 611\thinspace37 Brno, Czech
Republic\\\texttt{rosicky@math.muni.cz}\relax}

\begin{document}

\maketitle

\abstract\relax Are all subcategories of locally finitely presentable
categories that are closed under limits and $\lambda$-filtered colimits
also locally presentable? For full subcategories the answer is
affirmative. Makkai and Pitts proved that in the case $\lambda=\aleph_0$
the answer is affirmative also for all iso-full subcategories,
\emph{i.\thinspace e.}, those containing with every pair of objects all
isomorphisms between them. We discuss a possible generalization of this
from $\aleph_0$ to an arbitrary $\lambda$.\endabstract

\section{Introduction}

The lecture notes \cite{''4} of Horst Herrlich on reflections and coreflections in topology in 1968 with its introduction to categorical concepts
in Part II was the first text on category theory we read, and it has deeply influenced us. A quarter century later it was on Horst's impulse
that we solved the problem of finding full reflective subcategories of ${\mathbt{Top}}$ whose intersection is not reflective \cite{'1}. That paper 
has started years of intense cooperation of the authors. It is therefore with deep gratitude that we dedicate our paper to the memory of Horst Herrlich.
	
The above mentioned lecture notes are also the first reference for the fact that a full, reflective subcategory of a complete category is complete 
and closed under limits (9.1.2 in \cite{''4}; P.~J.~Freyd \cite{F} mentions this in Exercise 3F). In the present paper we study subcategories that 
are not necessarily full. They do not, in general, inherit completeness, see Example 2.1 below. However, if we restrict to complete subcategories, 
a necessary condition for reflectivity is that it be closed under limits. Now for full subcategories ${\mathscr K}$ of a locally presentable category 
an "almost" necessary condition is that ${\mathscr K}$ be closed under $\lambda$-filtered colimits for some $\lambda$, see Remark 2.3. Moreover, for full subcategories the Reflection Theorem of \cite1 (2.48) states the converse: a full subcategory is reflective if it is closed under limits and $\lambda$-filtered colimits. For non-full reflective subcategories the converse fails even in ${\mathbt{Set}}$: a subcategory of ${\mathbt{Set}}$ closed under limits and filtered colimits need not be reflective, see Example~\ref{E-2.2} below. A beautiful result was proved by Makkai and Pitts\thinspace\cite5 about \emph{iso-full} subcategories, \emph{i.\thinspace e.}, those containing every isomorphism of ${\mathscr L}$ with domain and codomain in the subcategory:

\subsection*{Theorem}\begingroup\em\csname@credit\endcsname[Makkai and
Pitts] Every iso-full subcategory of a locally finitely presentable
category closed under limits and filtered colimits ($\lambda=\aleph_0$)
is reflective.\endgroup\subsection*{}\indent

What about closedness under $\lambda$-filtered colimits for uncountable
$\lambda$? As an example, take ${\mathscr L}$ to be the category of
posets. Its subcategory of boolean $\sigma$-algebras and
$\sigma$-homomorphisms is iso-full and closed under limits and
$\aleph_1$-filtered colimits. Is it reflective? Yes, one can prove this
using Freyd's Special Adjoint Functor Theorem. However, we have not (in spite of
quite some effort) been able to answer the following general

\subsection*{Open problem} Is every iso-full subcategory of a locally
$\lambda$-presentable category that is closed under limits and
$\lambda$-filtered colimits reflective?\subsection*{}\indent

A substantial step in the proof of the above theorem due to Makkai and
Pitts was proving that the given subcategory is closed under elementary
subobjects. The main result of our paper is an ``approximation'' of the
affirmative answer to the above open problem based on introducing the
concept of a $\kappa$-elementary subobject (see Section~\ref{S-3}).
Every isomorphism is $\kappa$-elementary for all cardinals $\kappa$. And
conversely, every monomorphism that is $\kappa$-elementary for all
$\kappa$ is an isomorphism. Our ``approximate answer'' to the above open
problem substitutes for iso-fulness closedness under $\kappa$-elementary
subobjects:

\subsection*{Theorem}\begingroup\em Every subcategory ${\mathscr K}$ of
a locally $\lambda$-presentable category closed under limits,
$\lambda$-filtered colimits and $\kappa$-elementary subobjects for some
$\kappa\ge\lambda$ is reflective. And ${\mathscr K}$ is itself locally
$\lambda$-presentable.\endgroup\subsection*{}\indent

Another approximation concerns abelian categories ${\mathscr L}$: the
above open problem has an affirmative answer whenever the subcategory
contains all zero morphisms. Indeed, in this case the subcategory will
be proved to be full, thus our Reflection Theorem applies.

\subsection*{Acknowledgement} We are grateful to Saharon Shelah for a
consultation concerning a generalization of $\lambda$-good ultrafilters.
We hoped to use such generalized ultrafilters for a direct extension of
the proof of Makkai and Pitts from $\aleph_0$ to $\lambda$. However, the
result of the consultation is that the desired ultrafilters do not
exist. Pity!\subsection*{}\indent

\section{Reflective Subcategories}\label{S-2}

Recall that a subcategory ${\mathscr K}$ of a category ${\mathscr L}$ is
called reflective if the embedding $E:{\mathscr
K}\hookrightarrow{\mathscr L}$ has a left adjoint. The left adjoint
$F:{\mathscr L}\to{\mathscr K}$ is called a reflector, and the unit of
adjunction has components $$r_L:L\to FL\qquad\hbox{for $L\in{\mathscr
L}$}$$ called the reflections.

The subcategory ${\mathscr K}$ is said to be \emph{closed under limits}
in~${\mathscr L}$ if it has limits and the embedding $E$ preserves them.
(There is a stronger concept demanding that every limit in ${\mathscr
L}$ of a diagram in ${\mathscr K}$ lies in ${\mathscr K}$. For iso-full,
replete subcategories these two concepts coincide.) Analogously for
closure under finite limits, filtered colimits, \emph{etc.}

\subsection{Example} {\em An incomplete reflective subcategory of
${\mathbt{Pos}}$.} Let ${\mathscr K}$ have as objects all posets with a
least element~$0$ and a greatest element~$1$ such that $0\neq1$.
Morphisms are monotone functions preserving $0$ and~$1$. This category
is \begingroup

\hangindent1.666667\parindent\noindent\hskip\hangindent\llap{(i)\enspace
}iso-full,

\hangindent1.666667\parindent\noindent\hskip\hangindent\llap
{(ii)\enspace}reflective,

\noindent yet

\hangindent1.666667\parindent\noindent\hskip\hangindent\llap
{(iii)\enspace}incomplete: it does not have a terminal object.

\endgroup\hskip-\parindent\relax Indeed, an isomorphism clearly
preserves $0$ and $1$. A reflection of a poset is its embedding into the
poset with a new $0$ and a new $1$. For every object $K$ of ${\mathscr
K}$ more than one morphism leads from the $3$-chain to $K$, hence $K$ is
not terminal.

\subsection{Examples}\label{E-2.2} (1) {\em A subcategory ${\mathscr K}$ of
${\mathbt{Set}}$ that is \begingroup

\hangindent1.666667\parindent\noindent\hskip\hangindent\llap{(i)\enspace
}closed under limits and filtered colimits,

\noindent yet

\hangindent1.666667\parindent\noindent\hskip\hangindent\llap
{(ii)\enspace}not reflective.

\endgroup}\hskip-\parindent\relax Its objects are all sets
$X\times\{f\}$, where $X$ is a set and $f:{\mathrm{Ord}}\to X$ is a
function that is, from some ordinal onwards, constant. Morphisms from
$X\times\{f\}$ to $X'\times\{f'\}$ are all functions $h=h_0\times h_1$
where $h_0$ makes the triangle $$\xymatrix@!=12pt@C1pc{&{\mathrm{Ord}}
\ar[dl]_f\ar[dr]^{f'}&\\X\ar[rr]_{h_0}&&X'}$$ commutative.

The category ${\mathscr K}$ is equivalent to the category ${\mathscr
K}'$ of algebras of nullary operations indexed by ${\mathrm{Ord}}$ that
are from some ordinal onwards equal. The forgetful functor of ${\mathscr
K}'$ preserves limits and filtered colimits, whence it easily follows
that ${\mathscr K}$ is closed under limits and filtered colimits in
${\mathbt{Set}}$.

However, ${\mathscr K}$ does not have an initial object
(\emph{i.\thinspace e.}, the empty set has no reflection). Indeed, for
every object $K$ of ${\mathscr K}$ there exists an object $K'$ of
${\mathscr K}$ with ${\mathscr K}(K,K')=\emptyset$: for $K=X\times\{f\}$
with $f$ constant from $i$ onwards choose any $K'=X'\times\{f'\}$ with
$f'(i)\neq f'(i+1)$.

(2) {\em A subcategory ${\mathscr K}$ of
${\mathbt{Set}}$ that is \begingroup

\hangindent1.666667\parindent\noindent\hskip\hangindent\llap{(i)\enspace
}reflective and closed under filtered colimits,

\noindent yet

\hangindent1.666667\parindent\noindent\hskip\hangindent\llap
{(ii)\enspace}not iso-full.

\endgroup}\hskip-\parindent\relax Its objects are all sets
$X\times\{f\}$, where $X$ is a set and $f:X\to X$ is a function. Morphisms from
$X\times\{f\}$ to $X'\times\{f'\}$ are all functions $h$ such that $hf=f'h$.
 
The category ${\mathscr K}$ is equivalent to the category ${\mathscr K}'$ of algebras with one unary operation. The forgetful functor 
of ${\mathscr K}'$ is a right adjoint and preserves limits and filtered colimits, whence it easily follows
that ${\mathscr K}$ is reflective and closed under filtered colimits in ${\mathbt{Set}}$. But ${\mathscr K}$ is not iso-full.

\subsection{Definition} A non-full subcategory ${\mathscr K}$ of
${\mathscr L}$ is said to be \emph{closed under split subobjects} if for
every object $K$ of ${\mathscr K}$ and every split subobject of $K$
in~${\mathscr L}$ there exists a split monomorphism of ${\mathscr K}$
representing the same subobject.

Analogously closure under other types of subobjects is defined.

\proposition\label{P-2.4} An iso-full reflective subcategory is closed
under split subobjects iff it is full. \endproposition

\proof Let ${\mathscr K}$ be reflective in~${\mathscr L}$.

\smallskip(1)\enspace If ${\mathscr K}$ is full, then for every split
subobject $$m:L\to K,\quad e:K\to L,\quad em={\mathrm{id}},$$
in~${\mathscr L}$ with $K\in{\mathscr K}$ there exists a unique morphism
$\bar m$ of ${\mathscr K}$ such that the triangle
$$\xymatrix@!=12pt@C2.5pc@R2.5pc{L\ar[r]<2pt>^m\ar[d]_{r_L}&K\ar
@{-->}[l]<2pt>^e\\FL\ar[ru]_{\bar m}}$$ commutes. Then $r_L$ is inverse
to $e\bar m$. Indeed $$(e\bar m)\cdot r_L=em={\mathrm{id}}_L,$$ and for
the other identity use the universal property: from the equality
$$r_L\cdot(e\bar m)\cdot r_L=r_L={\mathrm{id}}_{FL}\cdot r_L$$ derive
$r_L\cdot(e\bar m)={\mathrm{id}}_L$.

Thus, $\bar m$ is the desired split monomorphism in ${\mathscr K}$: it
represents the same subobject as does $m$.

\smallskip(2)\enspace Let ${\mathscr K}$ be iso-full, reflective, and
closed under split subobjects. For every object $K$ of ${\mathscr K}$
the reflection of $K$ splits: there exists a unique morphism $e_K:FK\to
K$ of ${\mathscr K}$ with $e_K\cdot r_K={\mathrm{id}}_K$. Therefore,
$r_K$ represents the same subobject as does some split monomorphism
$\hat r_K:\hat K\to FK$ in ${\mathscr K}$; thus, we have an isomorphism
$i_K:\hat K\to K$ (in~${\mathscr L}$, thus in~${\mathscr K}$) for which
the triangle $$\xymatrix@!=12pt@C1pc{\hat K\ar[rr]^{i_K}\ar[dr]_{\hat
r_K}&&K\ar[dl]^{r_K}\\&FK}$$ commutes.

We prove that every morphism $f:K_1\to K_2$ of ${\mathscr L}$ between
objects of ${\mathscr K}$ lies in~${\mathscr K}$. The naturality square
$$\xymatrix@!=12pt@C2.5pc@R2.5pc{\hat K_1\ar@{-->}[r]^{i_{K_1}}\ar
@{-->}[dr]_{\hat r_{K_1}}&K_1\ar[r]^f\ar[d]_(.333333){r_{K_1}}&K_2\rlap
.\ar[d]<-2pt>_{r_{K_2}}\\&FK_1\ar[r]_{Ff}&FK_2\ar@{-->}[u]<-2pt>_e}$$
yields $$f=e\cdot Ff\cdot r_{K_1}=e\cdot Ff\cdot\hat
r_{K_1}\cdot{i_{K_1}}^{-1}.$$ This is a composite of morphisms of
${\mathscr K}$: as for ${i_{K_1}}^{-1}$ recall the iso-fulness.
\endproof

\subsection{Example} \emph{A non-full, iso-full reflective subcategory
of ${\mathbt{Pos}}$ that with every poset $K$ contains all split
subposets $K'$ as objects (but not necessarily the split monomorphism
$m:K'\to K$).}

Let ${\mathscr K}$ be the subcategory of all join semilattices with~$0$
and all functions preserving finite joins. This is clearly a non-full
but iso-full subcategory. A reflection of a poset $L$ is its embedding
$L\hookrightarrow{\mathrm{Id}}(L)$ into the poset of all ideals of $L$,
\emph{i.\thinspace e.}, down-closed and up-directed subsets, ordered by
inclusion.

Every split subobject (in ${\mathbt{Pos}}$) of a semilattice $K$,
$$m:K'\to K,\quad e:K\to K', \quad em={\mathrm{id}},$$ is itself a
semilattice. Indeed, the join of a finite set $M\subseteq K$ is easily
seem to be $e(\bigvee m[M])$. However, $m$ itself need not preserve
finite joins, as the following example demonstrates: $$
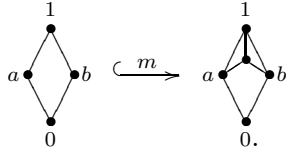
\xymatrix{{\xy*+{\xy/r1.5em/:(0,1)*\dir{*}="t";(+.5,0)*\dir{*}="1"**\dir
{-};(0,-1)*\dir{*}="b"**\dir{-};(-.5,0)*\dir{*}="0"**\dir{-};"t"**\dir
{-}\ar@{}"t"-/r1pt/;"t"+/r1pt/^{1\mathstrut}|{\dir{*}}\ar
@{}"b"-/l1pt/;"b"+/l1pt/^{0\mathstrut}|{\dir{*}}\ar
@{}"0"-/u1pt/;"0"+/u1pt/^{a\mathstrut}|{\dir{*}}\ar
@{}"1"-/d1pt/;"1"+/d1pt/^{b\mathstrut}|{\dir{*}}\POS\endxy}\endxy}\ar@{{
^(}->}[r]^m&{\xy*+{\xy/r1.5em/:(0,1)*\dir{*}="t";(+.5,0)*\dir
{*}="1"**\dir{-};(0,-1)*\dir{*}="b"**\dir{-};(-.5,0)*\dir{*}="0"**\dir
{-};"t"**\dir{-}\ar@{}"t"-/r1pt/;"t"+/r1pt/^{1\mathstrut}|{\dir{*}}\ar
@{}"b"-/l1pt/;"b"+/l1pt/^{0\rlap.\mathstrut}|{\dir{*}}\ar
@{}"0"-/u1pt/;"0"+/u1pt/^{a\mathstrut}|{\dir{*}}\ar
@{}"1"-/d1pt/;"1"+/d1pt/^{b\mathstrut}|{\dir{*}}\POS(0,.333333)*\dir
{*};"t"**\dir{-},"0"**\dir{-},"1"**\dir{-}\endxy}\endxy}}$$

\subsection{Remark} The reflection of an initial object of ${\mathscr
L}$ is clearly initial in~${\mathscr K}$. We say that ${\mathscr K}$ is
\emph{closed under initial objects} if ${\mathscr K}$ contains some
initial object of ${\mathscr L}$ that is initial in~${\mathscr K}$.

\proposition Let ${\mathscr L}$ be an abelian category. Every iso-full
subcategory closed under finite limits and initial objects is full.
\endproposition

\proof Let ${\mathscr K}$ be closed under finite limits and initial
objects. By Proposition\protect\unhbox\csname voidb@x\endcsname\penalty
\csname@M\endcsname\ \ref{P-2.4} we only need to prove that ${\mathscr
K}$ is closed under split subobjects. First notice that for every object
$K$ of ${\mathscr K}$ the coproduct injections of $K\oplus K$ lie
in~${\mathscr K}$. For example, the first coproduct injection $j:K\to
K\oplus K$ is the equalizer of the second product projection
$\pi_2:K\oplus K\to K$ and the zero morphism~$0$. Now $\pi_2$ lies
in~${\mathscr K}$ since ${\mathscr K}$ is closed under finite products,
and $0$ lies in~${\mathscr K}$ because ${\mathscr K}$ is also closed
under initial and terminal objects. Consequently, $j$ lies in~${\mathscr K}$.

For every split monomorphism $m:L\to K$, $K\in{\mathscr K}$, we prove
that $m$ lies in~${\mathscr K}$. There exists an object $B$ of
${\mathscr L}$ such that $m$ is the first coproduct injection of
$K=L\oplus B$. For the first coproduct injection $j:K\to K\oplus K$ it
follows that $m$ is the equalizer $$\xymatrix@!=12pt@C2.5pc@R2.5pc
{L\mathstrut\ar[r]^-m&K\mathstrut\ar[r]<2pt>^-j\ar
[r]<-2pt>_-{sj}&K\oplus K\mathstrut\rlap,}$$ where $$s:L\oplus B\oplus
L\oplus B\to L\oplus B\oplus L\oplus B$$ leaves the $L$-components
unchanged and swaps the $B$-components; shortly:
$$s=\langle\pi_1,\pi_4,\pi_3,\pi_2\rangle.$$ This follows easily from
$m$ being the first coproduct injection of $K=L\oplus B$. The morphism
$s$ lies in~${\mathscr K}$ since ${\mathscr K}$ is iso-full.
Consequently, $m$ lies in~${\mathscr K}$, as required. \endproof

\subsection{Example} {\em A non-full, reflective subcategory of
${\mathbt{Ab}}$. \/\rm It} has \begingroup

\hangindent1.666667\parindent\noindent\hskip\hangindent\llap{\enspace
}objects: all powers of the group $\ZZ$,

\hangindent1.666667\parindent\noindent\hskip\hangindent\llap{\enspace
}morphisms: all $\ZZ^u:\ZZ^I\to\ZZ^J$, where $u:J\to I$ is a function
and $\ZZ^u(h)=h\cdot u$.

\endgroup\hskip-\parindent\relax This subcategory ${\mathscr K}$ is not
full, since ${\mathscr K}(\ZZ^1,\ZZ^1)=\{{\mathrm{id}}_{\ZZ^1}\}$.
${\mathscr K}$ is reflective: the reflection of a group $G$ is its
canonical morphism $r:G\to\ZZ^{{\mathbt{Ab}}(G,\ZZ)}$, $r(x)(h)=h(x)$.
Indeed, for every group homomorphism $f:G\to\ZZ^I$ there exists a unique
function $u:I\to{\mathbt{Ab}}(G,\ZZ)$ with $f=\ZZ^u\cdot r$: put
$u(i)(x)=f(x)(i)$.

Observe that ${\mathscr K}$ is iso-full and closed under limits in $\mathbt{Ab}$
but it is not closed under initial objects. In fact, the initial object in ${\mathscr K}$
is $\ZZ^1$ and not $\ZZ^\emptyset$.

\subsection*{}\indent

We have asked, for a given locally $\lambda$-presentable category,
whether the iso-full subcategories closed under limits and
$\lambda$-filtered colimits are reflective in~${\mathscr L}$. Instead,
we can ask whether those subcategories are themselves locally
$\lambda$-presentable. This is an equivalent question:

\proposition\label{P-2.9} Let ${\mathscr L}$ be a locally $\lambda$-presentable
category. For iso-full subcategories ${\mathscr K}$ closed under limits
and $\lambda$-filtered colimits the following statements are equivalent:
\begingroup

\hangindent1.666667\parindent\noindent\hskip\hangindent\llap{(i)\enspace
}${\mathscr K}$ is reflective in~${\mathscr L}$

\noindent and

\hangindent1.666667\parindent\noindent\hskip\hangindent\llap
{(ii)\enspace}${\mathscr K}$ is a locally $\lambda$-presentable
category.

\endgroup\endproposition

\proof ii$\to$i. We apply the Adjoint Functor Theorem of the following
form proved in \cite1, Theorem~1.66: a functor between locally
presentable categories is a right adjoint iff it preserves limits and
$\lambda$-filtered colimits for some infinite cardinal~$\lambda$. We
conclude that the embedding $E:{\mathscr K}\to{\mathscr L}$ has a left
adjoint.

i$\to$ii. The left adjoint $F:{\mathscr L}\to{\mathscr K}$ of the
inclusion $E:{\mathscr K}\to{\mathscr L}$ preserves
$\lambda$-presentable objects. Indeed, given $L$ $\lambda$-presentable
in~${\mathscr L}$, we are to prove that $FL$ is $\lambda$-presentable
in~${\mathscr K}$. That is, ${\mathscr K}(FL,{\mskip1.5mu
\char"7B\mskip1.5mu })$ preserves colimits of $\lambda$-filtered
diagrams $D:{\mathscr D}\to{\mathscr K}$. For a colimit of $D$, its
image under $E$ is a colimit of $ED$ in~${\mathscr L}$. From the fact
that ${\mathscr L}(L,{\mskip1.5mu \char"7B\mskip1.5mu })$ preserves this
colimit it easily follows that ${\mathscr K}(FL,{\mskip1.5mu
\char"7B\mskip1.5mu })$ preserves the colimit of $D$ (using $F\dashv
E$).

The category ${\mathscr K}$ is cocomplete since ${\mathscr L}$ is, and
we prove that the objects $FL$, where $L$ ranges over
$\lambda$-presentable objects of ${\mathscr L}$, form a strong generator
of ${\mathscr K}$. Since these objects are $\lambda$-presentable
in~${\mathscr K}$ (and form a set up to isomorphism), it follows that
${\mathscr K}$ is locally $\lambda$-presentable by Theorem~1.20
of~\cite1. Thus, our task is to prove that for every proper subobject
$m:K\to K'$ in~${\mathscr K}$ there exists a morphism from some $FL$ to
$K'$, where $L$ is $\lambda$-presentable in~${\mathscr L}$, not
factorizing through $m$. We know that $m$ is not an isomorphism in
${\mathscr L}$ (since, ${\mathscr K}$ being iso-full, it would lie in
${\mathscr K}$). Since ${\mathscr L}$ is locally $\lambda$-presentable,
there exists a morphism $p:L\to K'$, $L$ $\lambda$-presentable in
${\mathscr L}$, that does not factorize through $m$. The unique morphism
$\bar p:FL\to K'$ of ${\mathscr K}$ with $p=\bar pr$, where $r$ is a
reflection, ^^ba^^ba does not factorize through $m$ either: given $u$
with $\bar p=mu$, we would have $p=mur$, a contradiction. \endproof

\section{Elementary Subobjects}\label{S-3}

We have mentioned the result of Makkai and Pitts that every iso-full
subcategory ${\mathscr K}$ of a locally finitely presentable category
${\mathscr L}$ closed under limits and filtered colimits is reflective.
A substantial step in the proof was to verify that in case ${\mathscr
L}$ is the category $\mathop{\mathbt{Str}}{\mit\Sigma}$ of structures of
some finitary, many-sorted signature ${\mit\Sigma}$, the given
subcategory ${\mathscr K}$ is closed under elementary subobjects. Recall
that a monomorphism $m:L\to K$ in $\mathop{\mathbt{Str}}{\mit\Sigma}$ is
called an \emph{elementary embedding} provided that for every formula
$\varphi(x_i)$ of first-order (finitary) logic with free variables $x_i$
and every interpretation $p(x_i)$ of the variables in~$L$ the following
holds. \begin{equation}\label{e-3.1} \hbox{$L$ satisfies
$\varphi(p(x_i))$ iff $K$ satisfies $\varphi(m(p(x_i)))$}.
\end{equation} (For many-sorted structures the variables also are
assigned sorts and interpretations are required to preserve sorts.)

We now consider the infinitary first-order logic ${\mathrm
L}_{\kappa\kappa}$, which allows conjunctions of fewer than $\kappa$
formulas and quantification over fewer than $\kappa$ variables. A
monomorphism $m:L\to K$ is called a \emph{$\kappa$-elementary} embedding
if $(\ref{e-3.1})$ holds for all formulas $\varphi(x_i)$ of ${\mathrm
L}_{\kappa\kappa}$.

\subsection{Example}\label{E-3.1} \smallskip(1)\enspace Every isomorphism is
$\kappa$-elementary for all cardinals $\kappa$.

\smallskip(2)\enspace The category of directed graphs is
$\mathop{\mathbt{Str}}{\mit\Sigma}$, where ${\mit\Sigma}$ consists of
one binary relation $R$. If $m:L\to K$ is a $\kappa$-elementary
embedding and $L$ has fewer than $\kappa$ vertices, then $m$ is an
isomorphism. Indeed, we can use the vertices of $L$ as variables $x_i$
($i\in I$); let $E$ be the set of all edges. The following formula
describes $L$: $$ \mathop{\,\bigvee\,}\limits_{\hidewidth(x_i,x_j)\in
E\hidewidth}R(x_i,x_j) \;\wedge\;\mathop{\,\bigvee\,}\limits_{\hidewidth
(x_i,x_j)\not\in E\hidewidth}\mathop\neg R(x_i,x_j)
\;\wedge\;\mathop{\mathop{\,\bigvee\,}\limits_{\hidewidth i,j\in
I\hidewidth}}\limits_{\hidewidth i\neq j\hidewidth}\mathop\neg x_i=x_j
\;\wedge\;\mathop{(\forall x)}\mathop{\,\bigvee\,}\limits_{\hidewidth
i\in I\hidewidth}x=x_i .$$ Since the formula holds in~$L$ for the
identity interpretation, it holds in~$K$ for $x_i\mapsto m(x_i)$. This
shows that $m$ is invertible.

Consequently, the only monomorphisms that are $\kappa$-elementary
embeddings for all $\kappa$ are isomorphisms.

\smallskip(3)\enspace More generally for every signature ${\mit\Sigma}$:
a morphism of $\mathop{\mathbt{Str}}{\mit\Sigma}$ is a
$\kappa$-elementary embedding for all $\kappa$ iff it is an isomorphism.

\subsection{Notation}\label{N-3.2} Recall that every locally
$\lambda$-presentable category ${\mathscr L}$ has a small full
subcategory ${\mathscr L}_\lambda$ representing all
$\lambda$-presentable objects up to isomorphism.

\smallskip(a)\enspace We denote by $${\mit\Sigma}_{\mathscr L}$$ the
many-sorted signature of unary operation symbols with

\hangindent1.666667\parindent\noindent\hskip\hangindent\llap{\enspace
}sorts = objects of ${\mathscr L}_\lambda$,

\noindent and

\hangindent1.666667\parindent\noindent\hskip\hangindent\llap{\enspace
}operation symbols = morphisms of ${\mathscr L}_\lambda^{\rm op}$.

\hskip-\parindent\relax Thus a morphism $f:s\to t$ of ${\mathscr
L}_\lambda^{\rm op}$ is a unary operation of input sort~$s$ and output
sort~$t$.

\smallskip(b)\enspace Every object $L$ of ${\mathscr L}$ defines a
${\mit\Sigma}_{\mathscr L}$-algebra $EL$: The underlying many-sorted set has components 
$$(EL)_s={\mathscr L}(s,L)\qquad\hbox{for all $s\in{\mathscr
L}_\lambda$}.$$ For an operation symbol $f:s\to t$ (morphism in
${\mathscr L}_\lambda^{\rm op}$) we define the operation of $EL$ by
precomposition with $f$: $$({\mskip1.5mu \char"7B\mskip1.5mu })\cdot
f:{\mathscr L}(s,L)\to{\mathscr L}(t,L).$$

\smallskip(c)\enspace Every morphism $h:L\to L'$ of ${\mathscr L}$
defines a homomorphism $Eh:EL\to EL'$ of ${\mit\Sigma}_{\mathscr
L}$-algebras: its components $(Eh)_s:{\mathscr L}(s,L)\to{\mathscr
L}(s,L')$ are given by postcomposition with $h$: $$h\cdot({\mskip1.5mu
\char"7B\mskip1.5mu }):{\mathscr L}(s,L)\to{\mathscr L}(s,L').$$

\lemma\label{L-3.3} \csname@credit\endcsname[see \cite{AR}, 1.26 and 1.27] For
every locally $\lambda$-presentable category ${\mathscr L}$ we have a
full embedding $$E:{\mathscr L}\to\mathop{\mathbt{Str}}{\mit\Sigma}_{
\mathscr L}$$ defined as above. The full subcategory $E({\mathscr L})$
is reflective and closed under $\lambda$-filtered colimits. \endlemma

\subsection{Definition}\label{D-3.4} A subcategory ${\mathscr K}$ of a
locally presentable category ${\mathscr L}$ is said to be \emph{closed
under $\kappa$-elementary embeddings} provided that there exists a
signature ${\mit\Sigma}$ and a full, reflective embedding $E:{\mathscr
L}\to\mathop{\mathbt{Str}}{\mit\Sigma}$ preserving $\kappa$-filtered colimits such that for every morphism
$m:L\to K$ of ${\mathscr L}$ with $K\in{\mathscr K}$ we have: if $Em$ is
a $\kappa$-elementary embedding, then $L$ and $m$ lie in~${\mathscr K}$.

\subsection{Remark} Any such a subcategory is iso-full and replete (see \ref{E-3.1}(1)). 

\section{Iso-Full Reflective Subcategories}

Recall from Makkai and Par\'e's \cite4 that a category ${\mathscr L}$ is
called $\lambda$-accessible if it has $\lambda$-filtered colimits and a
set of $\lambda$-presentable objects whose closure closure under
$\lambda$-filtered colimits is all of ${\mathscr L}$. An important
result of \cite4 is that for a signature ${\mit\Sigma}$ and a cardinal
$\lambda$ the category $\mathop{\mathbt{Elem}}_\lambda{\mit\Sigma}$ of
${\mit\Sigma}$-structures and $\lambda$-elementary embeddings is
$\kappa$-accessible for some $\kappa\ge\lambda$ (see also \cite{AR} 5.42).
Moreover, following the proof of \cite{AR} 5.42, a $\sigma$-structure
is $\kappa$-presentable iff its underlying set has cardinality $<\kappa$.

\theorem\label{T-4.1} Let ${\mathscr L}$ be a locally presentable
category. Every subcategory closed under limits, $\lambda$-filtered
colimits and $\lambda$-elementary embeddings for some $\lambda$ is
reflective in~${\mathscr L}$. \endtheorem

\proof Since instead of the given $\lambda$ every larger cardinal works
as well, we can assume without loss of generality that ${\mathscr L}$ is
locally $\lambda$-presentable and the embedding $E:{\mathscr
L}\to\mathop{\mathbt{Str}}{\mit\Sigma}$ of Definition~\ref{D-3.4}
preserves $\lambda$-filtered colimits.  

We are going to prove that the image $E({\mathscr K})$ is a reflective
subcategory of $\mathop{\mathbt{Str}}{\mit\Sigma}$. This implies that it
is reflective in $E({\mathscr L})$, and since $E$ defines an equivalence
of categories ${\mathscr L}$ and $E({\mathscr L})$, it follows that
${\mathscr K}$ is reflective in~${\mathscr L}$.

\smallskip(1)\enspace We first prove that every $\lambda$-presentable
${\mit\Sigma}$-structure $L\in{\mathscr L}$ has a reflection
in~${\mathscr K}$. By the preceding remark $\mathop{\mathbt{Elem}}_
\lambda{\mit\Sigma}$ is $\kappa$-accessible for some $\kappa\ge\lambda$, and
we take a set $\cal A$ of $\kappa$-presentable structures such that its
closure under $\kappa$-filtered colimits is all of that category.

We are going to prove that the slice category $L\mathbin/{\mathscr K}$
has an initial object (= reflection of $L$) by proving that the objects
$f:L\to K$ with $E(K)\in\cal A$ form a solution set. Since
$L\mathbin/{\mathscr K}$ is complete, Freyd's Adjoint Functor Theorem
then yields an initial object. Express $E(K)$ as a $k$-filtered colimit
$\bar k_i:\bar K_i\to E(K)$, $i\in I$, in $\mathop{\mathbt{Elem}}_\kappa
{\mit\Sigma}$ of objects $\bar K_i\in\cal A$. Since $\bar k_i$ is
$\kappa$-elementary and ${\mathscr K}$ is closed under
$\kappa$-elementary embeddings in~${\mathscr L}$ we obtain a
$\kappa$-filtered diagram $K_i$ in~${\mathscr K}$ whose image
$E(K_i)=\bar K_i$ is the given diagram and whose colimit cocone
$k_i:K_i\to K$ in~${\mathscr K}$ fulfils $Ek_i=\bar k_i$. Since this is
a $\kappa$-filtered colimit in~${\mathscr L}$ and $\kappa\ge\lambda$,
the morphism $f$ factorizes through some $k_i$:
$$\xymatrix@!=12pt@C2.5pc@R2.5pc{&K_i\ar[d]^{k_i}\\L\ar[ur]^{f'}\ar
[r]_f&K\rlap.}$$ The object $f':L\to K_i$ of $L\mathbin/{\mathscr K}$
lies in the specified set.

\smallskip(2)\enspace${\mathscr L}$ is reflective in~${\mathscr K}$.
Indeed, given an object $L$ of ${\mathscr L}$, express it as a
$\lambda$-filtered colimit $c_i:L_i\to L$ ($i\in I)$ of
$\lambda$-presentable objects $L_i$ of ${\mathscr L}$. By (1) we have
reflections $r_i:L_i\to FL_i$ in~${\mathscr K}$, and it is easy to see
that the objects $FL_i$ form a $\lambda$-filtered diagram in~${\mathscr
K}$ with a natural transformation having components $r_i$. Let
$r:L\to\mathop{\mathrm{colim}}_{i\in I} FL_i$ be a colimit of that
natural transformation. Then $r$ is a reflection of $L$ in~${\mathscr
K}$. Indeed, use closedness of ${\mathscr K}$ under $\lambda$-filtered
colimits: given a morphism $f:L\to K$ of ${\mathscr L}$ with
$K\in{\mathscr K}$, for every $i$ we get a unique morphism $f_i:FL_i\to
K$ of ${\mathscr K}$ with $fc_i=f_ir_i$: $$\xymatrix@!=12pt@C2.5pc
{L_i\ar[rr]^{r_i}\ar[dd]_{c_i}&&FL_i\ar[dl]_{f_i}\ar[dd]^{\bar
c_i}\\&K&\\L\ar[ur]^f\ar[rr]_r&&\mathop{\mathrm{colim}}FL_i\rlap.\ar
@{-->}[ul]_{\bar f}}$$

These morphisms form a cocone of the diagram of $FL_i$'s: If $d:L_i\to
L_j$ is a connecting morphism, the corresponding connecting morphism
$\bar d:FL_i\to FL_j$ is the unique morphism of ${\mathscr K}$ with
$\bar dr_i=r_jd$. Then in the following diagram
$$\xymatrix@!=12pt@C1pc{L_i\ar[rrrr]^d\ar[dr]^{r_i}\ar
[dd]_{c_i}&&&&L_j\ar[dl]_{r_j}\ar[dd]^{c_j}\\&FL_i\ar[rr]^{\bar
d}\ar[dr]_{f_i}&&FL_j\ar[dl]^{f_j}&\\L\ar[rr]_f&&K&&L\ar[ll]^f}$$ the
inner triangle commutes, as desired, because it lies in~${\mathscr K}$
and commutes when precomposed with the universal arrow $r_i$ (since
$c_i=c_jd$). Consequently, there exists a unique morphism $\bar
f:\mathop{\mathrm{colim}}FL_i\to K$ in~${\mathscr K}$ with $\bar
f\cdot\bar c_i=f_i$ (where the $\bar c_i$ form the colimit cocone). The
desired equality $$f=\bar fr:L\to K$$ follows from the fact that $c_i$
is a collectively epic cocone and the first diagram above commutes.

The uniqueness of $\bar f$ easily follows from the fact that the cocone
$\bar c_i$ is collectively epic in~${\mathscr K}$. \endproof

\subsection*{Problem.} Is any iso-full, reflective subcategory of a locally presentable category
which is closed under limits and $\lambda$-filtered colimits also closed under $\kappa$-elementary
embeddings for some $\kappa$?\subsection*{}\indent

Makkai and Pitts proved that this is true for $\lambda=\aleph_0$ (with $\kappa=\aleph_0$). 

\if01
Analogously to \ref{D-3.4} we say that a subcategory ${\mathscr K}$ of a
locally presentable category ${\mathscr L}$ is \emph{closed
under $\kappa$-elementary equivalence} provided that there exists a
signature ${\mit\Sigma}$ and a full, reflective embedding $E:{\mathscr
L}\to\mathop{\mathbt{Str}}{\mit\Sigma}$ preserving $\kappa$-filtered colimits such that for every objects
$L$ of ${\mathscr L}$ and $K\in{\mathscr K}$ we have: if $EL$ is
a $\kappa$-elementary equivalent to $EK$, then $L$ lies in~${\mathscr K}$.

Let us add that two $\Sigma$-structures are $\kappa$-elementary equivalent provided that they satisfy the same sentences of $L_{\kappa\kappa}$.

\proposition
Let $\mathscr K$ be a replete and iso-full, reflective subcategory of a locally presentable category $\mathscr L$ closed under limits and $\lambda$-filtered colimits for some regular cardinal $\lambda$. Then $\mathscr K$ is closed in $\mathscr L$ under $\kappa$-elementary equivalence for some $\kappa$.
\endproposition
\proof
Following \ref{P-2.9}, $\mathscr K$ is locally presentable. Let $E:\mathscr L\to\mathbt{Str}{\mit\Sigma}$ be a full, reflective embedding
and $E':\mathscr K\to\mathbt{Str}{\mit\Sigma}$ its restriction on $\mathscr K$. There is a regular cardinal $\mu\geq\lambda$ such that both
$\mathscr K$ and $\mathbt{Str}{\mit\Sigma}$ are locally $\mu$-presentable, $E'$ preserves $\mu$-presentable objects and $\mu$-presentable objects
in $\mathbt{Str}{\mit\Sigma}$ are precisely objects of cardinality $<\mu$ (see \cite{AR} 5.42 and \cite4 2.4.9). Let $\mathscr A=(\mathbt{Str}{\mit\Sigma})_\mu$ 
(see Notation \ref{N-3.2}). Since $E'$ preserves $\mu$-presentable objects, there is a subcategory $\mathscr B\subseteq\mathscr A$ which is equal
to $\mathscr K_\mu$. For any $\mu$-directed diagram $D:I\to\mathscr B$, denote $D_{ij}:i\to j$ its morphisms, and consider the formula $\varphi_D$
$$
\bigwedge_{i\in I}(\exists x_i)\pi^+_{Di}(x_i)\wedge \bigwedge_{i\leq j} (x_jD_{ij}=x_i)\wedge\bigwedge_{A\in\mathscr A}(\forall y)(\pi^+_A(y)\to
$$
$$
(\bigvee_{t:A\to Di} x_it=y\wedge\bigwedge_{(t_1,t_2)}(x_{i_1}t_1=y\wedge x_{i_2}t_2=y\to (\bigvee_{(s_1,s_2)}(x_is_1t_1=y)) .
$$
Here $\pi^+_A$ is the positive diagram of $A\in\mathscr A$ (see \cite{AR} 5.33), $(t_1,t_2)$ is a pair of morphisms $t_1:A\to Di_1$, $t_2:A\to Di_2$ and
$(s_1,s_2)$ is a pair of morphisms $s_1:Di_1\to Di$, $s_2:Di_2\to Di$ with $s_1t_1=s_2t_2$. Since $\mathscr K$ is replete, for a $\Sigma$-structure $L$ 
we have $L\models\varphi_D$ with $x_i$ being interpreted as $d_i$ if and only if $d_i:Di\to L$ is a cocone for $D$ whose colimit is $L$ because $d_i:Di\to L$ 
is cofinal in the canonical cocone $h:A\to L$, $A\in\mathscr A$ (with $y$ being interpreted as $h$). For any $K\in\mathscr K$, there is $D:I\to\mathscr B$ such that $K\models\varphi_D$. Since the quantifier rank of $\varphi$ is $2$, we have $L\models\varphi_D$ for any $L\equiv^2_{\infty\mu} K$ where 
$\equiv^2_{\infty\mu}$ denotes elementary equivalence with respect to formulas of the quantifier rank $\leq 2$ in the logic $L_{\infty\mu}$ allowing arbitrary conjuctions and quantification over fewer that $\mu$ variables (see \cite{D}, p. 319). Let $\kappa = 2^{(2^{\mu^-\cdot |\Sigma|})}$ where $(\alpha^+)^- = \alpha$ and $\beta^- =\beta$ if $\beta$ is not of the form $\alpha^+$. By Benda \cite{B} (see also \cite{D}, p. 344) $L\equiv^2_{\infty\mu} K$ if and only if 
$L\equiv^2_{\kappa\mu} K$. Thus any $L\equiv_{\kappa\mu} K$ belongs to $\mathscr K$. Consequently, $\mathscr K$ is closed in $\mathscr L$ under 
$\kappa$-elementary equivalence.
\endproof

\fi
\subsection{Example} We apply Theorem \ref{T-4.1} to Kan injectivity, a concept
in order-enriched categories ${\mathscr L}$ introduced by
Escardo\thinspace\cite{'4} for objects and by Carvalho and
Sousa\thinspace\cite3 for morphisms.

An object $K$ of ${\mathscr L}$ is said to be \emph{Kan-injective}
w.\thinspace r.\thinspace t.\ a morphism $h:X\to Y$ provided that every
morphism $f:X\to K$ has a left Kan extension $f/h$ along $h$. That is,
there is a morphism $f/h:Y\to K$ that fulfils
\begin{equation}\label{e-5.1} f\le(f/h)\cdot h\end{equation} and is
universal with this property, \emph{i.\thinspace e.}, for every morphism
$g:Y\to K$ \begin{equation}\label{e-5.2} f\le g\cdot
h\quad\hbox{implies}\quad f/h\le g.\end{equation} Carvalho and Sousa
introduced in \cite{3} the category ${\mathscr L}\mathop{\mathrm{Inj}}({
\mathscr H})$, for every class ${\mathscr H}$ of morphisms of ${\mathscr
L}$: Objects are the those objects of ${\mathscr L}$ that are
Kan-injective w.\thinspace r.\thinspace t.\ every $h\in{\mathscr H}$.
Morphisms $u:K\to K'$ are those morphisms of ${\mathscr L}$ that
\emph{preserve Kan extensions}: given $h:X\to Y$ in~${\mathscr H}$, we
have $$u\cdot(f/h)=(uf)/h\qquad\hbox{for all $f:X\to K$}.$$

For example, in ${\mathscr L}={\mathbt{Pos}}$ let $h:X\to Y$ be the
embedding of a two-element discrete poset $X$ into $Y=X\cup\{t\}$, $t$ a
top element. Then ${\mathscr L}\mathop{\mathrm{Inj}}\{h\}$ is the
subcategory of join semilattices and their homomorphisms.

\subsection{Remark} Recall that an \emph{order-enriched locally
$\lambda$-presentable} category is a category that is locally
$\lambda$-presentable and enriched over ${\mathbt{Pos}}$ in such a way
that a parallel pair $f_1,f_2:L\to L'$ fulfils $f_1\le f_2$ whenever for
every morphism $u:K\to L$, $K$ $\lambda$-presentable, we have $f_1\cdot
u\le f_2\cdot u$. \endgraf\subsection*{}\indent\endgraf The following
result is proved in~\cite2, but the present proof is simpler.

\proposition Let $\mathscr H$ be a set of morphisms of an order-enriched
locally presentable category. Then ${\mathscr L}\mathop{\mathrm{Inj}}({\mathscr H})$ is a
reflective subcategory.
\endproposition

\proof Since ${\mathscr H}$ is a set, there exists a cardinal $\lambda$
such that the given category ${\mathscr L}$ is locally
$\lambda$-presentable and domains and codomains of members of ${\mathscr
H}$ are all $\lambda$-presentable. As can be seen rather easily it
follows that ${\mathscr L}\mathop{\mathrm{Inj}}({\mathscr H})$ is closed
under $\lambda$-filtered colimits in~${\mathscr L}$, see \cite{ASV} for
details. A proof that ${\mathscr L}\mathop{\mathrm{Inj}}({\mathscr H})$
is closed under limits can be found in~\cite{CS}. Thus according to
Theorem\protect\unhbox\csname voidb@x\endcsname\penalty\csname
@M\endcsname\ \ref{T-4.1} it is sufficient to prove closedness under
$\lambda$-elementary embeddings.

Let $\bar{\mit\Sigma}_{\mathscr L}$ be the extension of the signature
${\mit\Sigma}_{\mathscr L}$ of Notation~\ref{N-3.2} by a binary relation
symbol $\le_s$ for every sort $s\in{\mathscr L}_\lambda$. And let $\bar
E:{\mathscr L}\to\mathop{\mathbt{Str}}\bar{\mit\Sigma}_{\mathscr L}$ be
the extension of the full embedding of Lemma\protect\unhbox\csname
voidb@x\endcsname\penalty\csname@M\endcsname\ \ref{L-3.3} by
interpreting, for every object $L$ in~${\mathscr L}$, the symbol $\le_s$
as the given partial order on ${\mathscr L}(s,L)$. For every member
$h:s\to t$ of ${\mathscr H}$, since $s$ and $t$ lie in~${\mathscr
L}_\lambda$, we have the unary operation symbol $h:t\to s$ in
$\bar{\mit\Sigma}_{\mathscr L}$ (interpreted as precomposition with
$h$). The following formula, with variables $x$ of sort $s$ and $y$, $z$
of sort $t$, expresses Kan-injectivity w.\thinspace r.\thinspace
t.\ $h$: $$\psi_h \mathrel{\;=\;} \mathop{(\forall x)} \mathop{(\exists
y)} \left([x\le_sh(y)]\mathbin{\,\wedge\,} \mathop{(\forall z)}
\left([x\le_sh(z)]\to[y\le_tz]\right)\right).$$ Indeed, an object $L$ of
${\mathscr L}$ has the property that $\psi_h$ holds in $\bar EL$ iff for
every element of sort $s$, \emph{i.\thinspace e.}, every morphism
$f:s\to L$ in~${\mathscr L}$, there exists an element of sort $t$,
\emph{i.\thinspace e.}, a morphism $\bar f:t\to L$, such that
(i)~$f\le\bar f\cdot h$ and (ii)~for every interpretation of $z$,
\emph{i.\thinspace e.}, every morphism $g:t\to L$, if $f\le g\cdot h$,
then $\bar f\le g$. This tells us precisely that $\bar f=f/h$.

Now let $m:L\to K$ be a morphism of ${\mathscr L}$ with $K$
Kan-injective w.\thinspace r.\thinspace t.\ ${\mathscr H}$ and with
$\bar Em$ a $\lambda$-elementary embedding. We prove that $L$ is also
Kan-injective, and that $m$ preserves Kan extensions. Given
$h\in{\mathscr H}$ we know that $\bar EK$ satisfies $\psi_h$, and this,
since $\psi_h$ has no free variables, implies that also $\bar EL$
satisfies $\psi_h$. That is, $L$ is Kan-injective w.\thinspace
r.\thinspace t.\ $h$ (for all $h\in{\mathscr H}$). Next consider the
formula $\psi_h'(x,y)$ with free variables $x$ and $y$ obtained by
deleting the two quantifiers at the beginning of $\psi_h$. An
interpretation of $\psi_h'(x,y)$ in $\bar EL$ is a pair of morphisms
$f:s\to L$ and $\bar f:t\to L$ with $\bar f=f/h$. Analogously for
interpretations in $\bar EK$. The assumption that $\bar Em$ is a
$\lambda$-elementary embedding tells us: for an interpretation $f,f/h$
in $\bar EL$, it follows that $m\cdot f,m\cdot(f/h)$ is an
interpretation in $\bar EK$, that is, $m\cdot(f/h)=(m\cdot f)/h$, as
required. \endproof

\begin{references*}

\bibitem{'1} J.~Ad\'amek and J.~Rosick\'y, Intersections of reflective subcategories,
Proc. Amer. Math. Soc. 103 (1988), 710-712.

\bibitem1\immediate\write\csname@auxout\endcsname{\string\bibcite
{AR}{\the\csname c@reflister\endcsname}} J.~Ad\'amek and J.~Rosick\'y,
\emph{Locally Presentable and Accessible Categories}, London
Mathematical Society Lecture Note Series 189, Cambridge University
Press 1994.

\bibitem2\immediate\write\csname@auxout\endcsname{\string\bibcite
{ASV}{\the\csname c@reflister\endcsname}} J.~Ad\'amek, L.~Sousa and
J.~Velebil, Small Kan-injectivity subcategories, manuscript.

\bibitem{B} M. Benda, Reduced products and non-standard logic, J. Symb. Logic 34 (1969), 424-436.

\bibitem{3}\immediate\write\csname@auxout\endcsname{\string\bibcite
{CS}{\the\csname c@reflister\endcsname}} M.~Carvalho and L.~Sousa,
Order-preserving reflectors and injectivity, \emph{Topology Appl.} 158
(17) (2011), 2408--2422.

\bibitem{D} M. A. Dickmann, \emph {Large Infinitary Languages}, North-Holland 1975.

\bibitem{'4} M.~H.~Escard\'o, Injective spaces via the filter monad,
Proc.\ 12th Summer Conference on General Topology and Applications
(1997), 97--110.

\bibitem{F} P.~J.~Freyd, \emph{Abelian Categories}, Harper and Row, New York 1964.

\bibitem{''4} H.~Herrlich, Topologische Reflexionen und Coreflexionen,
Lecture Notes in Math. 78, Springer-Verlag 1968.

\bibitem4 M.~Makkai and R.~Par\'e, \emph{Accessible Categories: The
Foundations of Categorical Model Theory}, Contemporary Mathematics 104,
American Mathematical Society 1989.

\bibitem5 M.~Makkai and A.~M.~Pitts, Some results on locally finitely
presentable categories, \emph{Trans.\ Amer.\ Math.\ Soc.} 299 (1987),
473--496.

\end{references*}

\end{document}